\documentclass[a4paper]{scrartcl}
\usepackage{amstext}
\usepackage{amsmath}
\usepackage{amssymb}
\usepackage{amsthm}
\usepackage{url}
\usepackage{todonotes}

\theoremstyle{plain}

\newtheorem{Thm}{Theorem}

\newtheorem{Cor}[Thm]{Corollary}

\newcommand{\NN}{\mathbb{N}}
\newcommand{\RR}{\mathbb{R}}

\newcommand{\eps}{\varepsilon}

\begin{document}
\title{Infinity in computable probability}
\subtitle{Logical proof that William Shakespeare probably was not a dactylographic monkey}
\author{Maarten McKubre-Jordens\footnote{Adjunct Fellow, School of Mathematics \& Statistics, University of Canterbury, New Zealand.}\hspace{2mm} \& Phillip L.\! Wilson\footnote{Corresponding Author: phillip.wilson@canterbury.ac.nz. School of Mathematics \& Statistics, University of Canterbury, New Zealand.}}
\maketitle

\section{Introduction}
Since at least the time of Aristotle \cite{aristotle}, the concept of combining a finite number of objects infinitely many times has been taken to imply certainty of construction of a particular object. In a frequently-encountered modern example of this argument, at least one of infinitely many monkeys, producing a character string equal in length to the collected works of Shakespeare by striking typewriter keys in a uniformly random manner, will with probability one reproduce the collected works. In the following, the term ``monkey'' can (naturally) refer to some (abstract) device capable of producing sequences of letters of arbitrary (fixed) length at a reasonable speed.

Recursive function theory is one possible model for computation; Russian recursive mathematics is a reasonable formalization of this theory \cite{varieties}.\footnote{The history of the relationship between classical logic and computation is long and complex, and beyond the scope of this paper. These results do indicate that mathematics using constructive logics, such as the Russian recursive mathematics used here, seems to be more suited to the simulation of the work of a computer than classical logic.} Here we show that, surprisingly, within recursive mathematics it is possible to assign to an infinite number of monkeys probabilities of reproducing Shakespeare's collected works in such a way that while it is impossible that no monkey reproduces the collected works, the probability of \emph{any} finite number of monkeys reproducing the works of Shakespeare is \emph{arbitrarily} small. The method of assigning probabilities depends only on the desired probability of success and not on the size of any finite subset of monkeys. 

Moreover, the result extends to reproducing all possible texts of any finite given length. However, in the context of implementing an experiment or simulation computationally (such as the small-scale example in \cite{paignton}; see also \cite{complete}), the fraction among all possible probability distributions of such \emph{pathological} distributions is vanishingly small provided sufficiently large samples are taken.

\section{The classical experiment}

The classical infinite monkey theorem \cite{borel,eddington} can be stated as follows: given an infinite amount of time, a monkey hitting keys on a typewriter with uniformly random probability will almost certainly type the collected works of William Shakespeare \cite{shakespeare}. We use a slightly altered (but equivalent) theorem involving an infinite collection of monkeys, and give an intuitive direct proof.

Let a string of characters of length $w \in \NN^+$ over a given alphabet $A$ (of size $|A|$, including punctuation) be called a \emph{$w$-string}. For example, ``\emph{banana}'' is a 6-string over the alphabet $\{a,b,n\}$. Suppose each monkey is given a computer keyboard
%\footnote{In this modern day and age, we won't skimp on what we give our monkey.}
 with $|A|$ keys, each corresponding to a different character. Suppose also that the experiment is so contrived that each monkey will type its $w$-string in finite time. 

\begin{Thm}\label{thm:imt}
At least one of infinitely many monkeys typing $w$-strings, as described in the previous paragraph, will almost certainly produce a perfect copy of a target $w$-string in finite time.
\end{Thm}

\begin{proof}
Recall that for this theorem, the probability that any given monkey strikes any particular key is uniformly distributed. Let the target $w$-string be $T_w$. %Fix a monkey, say Albert \cite{nasa}. Then 
The chance of a given monkey producing $T_w$ is simply the probability of him typing each character of the target text in the correct position, or
\[
\underset{w \text{ terms}}{\underbrace{\frac{1}{|A|} \times \frac{1}{|A|} \times \dots \times \frac{1}{|A|}}} = \left(\frac{1}{|A|}\right)^w\ .
\]
Therefore, the probability that a given monkey \emph{fails} to produce $T_w$ is
\[
1- \left(\frac{1}{|A|}\right)^w\ .
\]

Now, if we examine the output of $m$ monkeys, then the probability that \emph{none} of these monkeys produces $T_w$ is
\[
T_w(m) = \left(1- \left(\frac{1}{|A|}\right)^w\right)^m.
\]
Therefore, the probability that at least one monkey of $m$ produces $T_w$ is
\[
P(m) = 1-T_w(m)\ .
\]
Now
\[
\lim_{m \to \infty} P(m) = 1.
\]
\end{proof}
In other words, as the number of monkeys tends to infinity, at least one will almost certainly produce the required string.

However, any real-world experiment that attempts to show this will, unless the target $w$-string and $|A|$ are quite small, be very likely to fail, since the probabilities involved are tiny. For example, taking the English alphabet (together with punctuation and capitalization) to have 64 characters, a simple computation shows that, if the monkeys are typing 6-strings, the chances of a monkey typing ``banana'' correctly are
\begin{equation}\label{eq:tiny}
\left(\frac{1}{64}\right)^6 = \frac{1}{68719476736} \approx 1.5 \times 10^{-11}.
\end{equation}
If it takes one second to check a single monkey's output, then the number of seconds that will elapse before we have a 50\% chance of finding a monkey that has typed ``banana'' correctly is outside the precision of typical computing software. Of course, if some monkeys have a preference for typing a certain letter more often than others --- say `a' --- then this probability can be much larger. Indeed, it is non-uniformity among monkeys
%\footnote{And why \emph{should} we demand that all monkeys are created equal?}
 that we exploit to derive our main result in \S\ref{sec:recursive}.
 
Results such as (\ref{eq:tiny}) have been interpreted \cite[p.53]{kk} as saying that ``The probability of [reproducing the collected works of Shakespeare] is therefore zero in any operational sense\ldots''. In \S\ref{sec:recursive} we show that this probability can be made arbitrarily small in any sense, operational or otherwise.

\section{A simple, classical non-uniform version}\label{sec:class}

What if the monkeys do not necessarily strike their keys in a uniformly distributed manner? In this case, we might prescribe a certain probability for a particular monkey to type a particular $w$-string (and this probability need not be the same from one monkey to the next). Before we reach our main result, we outline a non-uniform classical probability distribution such that for any $\eps>0$ the probability of success by monkey $m$ is arbitrarily small, but with the probability of failure still zero. If we allow our probability distribution to be a function of $m$ as well as $\eps$ then the following distribution will suffice:
\[
1-p_k(m,\eps)=\delta(m-k)(\eps-\sigma)+\delta(m+1-k)(1-\eps+\sigma)\ ,
\] 
where $p_k$ is the probability of failing at monkey $k$, the Dirac delta function $\delta(s)=1$ for $s=0$ and is zero otherwise, and $0<\sigma<\eps$. Here, the probability of finding the target $w$-string at or before the $m^{\text{th}}$ monkey, $P(m)$, is less than the prescribed $\eps$, but success is still certain --- we need merely look at $m+1$ monkeys.

In the following section we show that, surprisingly, this can be achieved with a probability distribution dependent only on $\eps$, and not on $m$. That is, it is possible to produce a computable distribution so that, while each monkey produces Shakespeare's works with nonzero probability, actually finding the culprit among \emph{any} finite subcollection is very unlikely. To do so, we invoke a result from recursive mathematics.

\section{The successful monkey is arbitrarily elusive}\label{sec:recursive}

Within recursive mathematics, there is a theorem sometimes referred to as the \emph{singular covering theorem}, originally proved by Tseitin and Zaslavsky (1956), and independently by Kreisel and Lacombe (1957) (see \cite{kushner}): given a compact set $K$, for every positive $\varepsilon$, one can construct a computable open rational $\varepsilon$-bounded covering of $K$.\footnote{Related theorems with detailed proofs and discussion were published by Tseitin and Zaslavsky in \cite{zaslavsky}. We hasten to add that while this may seem esoteric, the results really provide commentary on much more mainstream ideas such as computer simulations, since constructive logics are much more suited to theorizing about these.} It can be restricted to the interval $[0,1]$ as follows:

\begin{Thm}\label{thm:sct}
For each $\eps > 0$ there exists a sequence $(I_k)_{k = 1}^\infty$ of bounded open rational intervals in $\RR$ such that
\begin{enumerate}
\item[\textrm{(i)}] $[0,1] \subset \bigcup_{k=1}^\infty I_k$, and
\item[\textrm{(ii)}] $\sum_{k=1}^n |I_k| < \eps$ for each $n \in \NN^+$.
\end{enumerate}
\end{Thm}
Our principal result, Theorem \ref{thm:eps}, follows from this theorem, and highlights the tension between classical probability theory and its constructive counterpart as outlined in \cite{chan}.

To set up our principal theorem, we first define $M$ to be an infinite, enumerable set of monkeys (the \emph{monkeyverse}), and for any natural number $m$ the \emph{$m$-troop} of monkeys to be the first $m$ monkeys in $M$.  Note that for any given monkey it is decidable whether that monkey has produced a given finite target string.

\begin{Thm} \label{thm:eps}
Given a finite target $w$-string $T_w$ and a positive real number $\eps$, there exists a computable probability distribution on $M$ of producing $w$-strings such that:
\begin{enumerate}
\item[\textrm{(i)}] the classical probability that no monkey in $M$ produces $T_w$ is $0$; and
\item[\textrm{(ii)}] the probability of a monkey in any $m$-troop producing $T_w$ is less than $\eps$.
\end{enumerate}
\end{Thm}
\begin{proof}
 Suppose that the hypotheses of the theorem are satisfied. As above, let $P(m)$ be the probability that a monkey in the $m$-troop has produced $T_w$, and let $p_k$ be the probability that the $k^{\text{th}}$ monkey has \emph{not} produced $T_w$. Then
\[
P(m) = 1-\prod_{k=1}^m p_k.
\]
Given $0 < \eps < 1$, compute $\eps_0 = -\log (1-\eps)$. For this $\eps_0$, construct the singular cover $(I_k)_{k=1}^\infty$ as per Theorem \ref{thm:sct}. %Some care must be taken: the singular cover does \emph{not} guarantee that $\sum_{k=1}^\infty |I_k| = 1$; this sum may in fact be bigger than one. To prevent this, take the \emph{restriction} of the singular cover to $[0,1]$. 
Then set
\[
p_k = \exp\left(-|I_k|\right).
\]
To prove (i), observe that $0<p_k<1$ for each $k$. The monotone convergence theorem now ensures that the product $\prod_{k=1}^m p_k$ \emph{classically} tends to $0$, hence it is (classically) impossible that no monkey produces $T_w$.

On the other hand, we have (computably)
\[
-\log(p_k) = |I_k|,
\]
whence, by the singular covering theorem,
\[
\sum_{k=1}^m -\log(p_k) = \sum_{k=1}^m |I_k| < \eps_0 = -\log(1-\eps)
\]
for all $m \in \NN^+$. Some rearranging shows that
\[
\log\left(\prod_{k=1}^m p_k \right) = \sum_{k=1}^{m} \log(p_k) > \log(1-\eps)
\]
and hence
\[
\prod_{k=1}^m p_k > 1- \eps.
\]
Then the probability of any member of the $m$-troop producing $T_w$ is
\[
P(m) = 1 - \prod_{k=1}^m p_k < \eps
\]
for \emph{any} positive natural number $m$. This proves (ii).
\end{proof}

Thus, the chances of us actually \emph{finding} the monkey that produces the collected works of Shakespeare can be made arbitrarily small, and the classical intuition that, since we have an infinite number of monkeys, Shakespeare's works must be typed by \emph{some} monkey is of no help in \emph{locating} the successful monkey.

%\section{What have we actually done?}

%In Theorem \ref{thm:eps} we have allowed that the probability of producing $T_w$ varies between monkeys.
%(given a $w$-string, the probability of Albert producing it could be different from the probability of William producing it). 
%This crucial difference from the uniform case allows us to assign probabilities to the monkeyverse of producing the target $w$-string in such a way that:
%\begin{itemize}
%\item the probability that \emph{no} monkey has produced the target $w$-string is 0; and
%\item the probability that we can \emph{find} the monkey that has produced the target $w$-string is arbitrarily small, no matter how big a troop of monkeys we search.
%\end{itemize}
We emphasize that, in contrast to the case in \S\ref{sec:class}, the pathological distribution in Theorem \ref{thm:eps} does \emph{not} depend on $m$, the size of the troop we search.\footnote{Contrasting the classical with the computational view in the same proof may prove counterintuitive. We are hoping to shed light on why the \emph{intuitive} result---that it is (in the classical abstract world) impossible that \emph{no} monkey produces Shakespeare's works---clashes with the fact that it may be incredibly difficult (in the concrete computational world) to nail the cheeky monkey that did it.

What sense to make of the product $\prod p_k$ of monkeys failing to produce Shakespeare classically tending to $0$? The problem here is the \emph{rate} at which it does so---this rate is computationally untractable.}

One might argue that it is easy to assign probabilities in such a way that any finite search will almost certainly not yield the monkey that produced it --- by letting each monkey produce the target $w$-string with probability zero. However, in this case, \emph{no} monkey will produce it. Our theorem shows that, even in the case where it is (classically) impossible that no monkey produces the target, it is still possible to make the probability of finding the monkey that accomplishes the necessary task arbitrarily small.

\section{Target-free writing}

%In terms of practical fiction writing, o
One criticism of the above line of reasoning is that the experimenter requires knowledge of the target. There, the output of each monkey was tested against the collected works of Shakespeare: only if every character matched would it pass the test. However, suppose now that 
%in some dystopian future where the works of Shakespeare have been lost, but monkeys remain plentiful\footnote{For some, especially monkeys, this scenario may in fact be \emph{utopian}.}, 
we wish to recreate Shakespeare's work armed only with knowledge of the total character length in some alphabet. That is, we know that we require one of the $|A|^w$ possible $w$-strings. Can we guarantee to complete the list (without repetition) and therefore recreate the collected works of Shakespeare (somewhere)? We note that the list can be shortened by checking for grammar etc.\footnote{Truncating the list in this way may be desirable in order to avoid being overwhelmed by ``meaningless cacophonies, verbal farragoes, and babblings'' \cite{borges}.}; here we consider the worst case of the complete list, without repetition, of $w$-strings.

\begin{Cor}
Any list of finite strings is completed in finite time with arbitrarily small probability.
\end{Cor}

The proof relies on applying Theorem \ref{thm:eps} multiple times using standard calculations.

\section{Pathological distributions are arbitrarily rare}

At first sight, Theorem \ref{thm:eps} might appear to destroy any hope of finding the successful monkey. However, we have the following:
\begin{Thm}
The probability that the probability distribution on the monkeyverse is constructed in such a way as to make the constructive probability of finding the desired monkey arbitrarily small, is arbitrarily small.
\end{Thm}
\begin{proof}
Given $0 < \eps < 1$, in order for the probability distribution to be pathological, the probability of \emph{any} monkey in the $m$-troop outputting $T_w$ cannot exceed $\eps$. Therefore the fraction of pathological distributions over an $m$-troop is at most $\eps^m$, and \[\lim_{m \to \infty} \eps^m = 0.\]
\end{proof}
In short, we can make the fraction of pathological distributions arbitrarily small if we search sufficiently large $m$-troops. Here, then, is an a priori justification for large sample sizes in the case of computational simulations.
%\begin{proof}
%With the notation from Theorem \ref{thm:eps}, if for the $k^{\text{th}}$ monkey $\log (p_k)\geqslant 1-\eps$ then $1-\prod_{k=1}^mp_k>\eps$ for any $m\geqslant k$. Consequently the fraction of non-pathological probability distributions over an $m$-troop is \[F_m\geq 1-\left(\frac{\eps}{1-\eps}\right)^m\ .\] If $\eps<0.5$, \[\lim_{m\rightarrow\infty}F_m=1\ .\]
%\end{proof}
%In other words, from a frequentist point of view, the probability of encountering a pathological distribution in a computer-implemented simulation is vanishingly small.

\section{Discussion and further work}

Recall that, throughout this paper, we take the term ``monkey'' to refer to some device capable of producing arbitrary but finite sequences of letters --- computers satisfy this criterion. The theorems presented in this paper therefore have implications for computer simulations. In particular, when performing simulations of a probabilistic nature, the experimenter needs to ensure that pathological distributions do not arise, or arise rarely enough to provide a measure of confidence in the conclusion.

%because it is related to such deep facts as that from the classical viewpoint the computable real numbers have measure 0 (and all finite texts “printed by monkeys” correspond to rational numbers/ subsets of computable real numbers). There is also a philosophical issue of intuitionistic free choice sequences as opposed to computable sequences, of various versions of König’s lemma (constructively valid and not) etc. I do not say that the authors in this paper should go deep into discussion of these issues, but an outline of the context would strengthen the point raised in their paper beyond a simple curiosity. 

It should also be noted that the classical non-uniform distribution presented above suggests that a pathological situation can never be ruled out with certainty, since if the experimenter tests just one more monkey, the result may be vastly different than observed earlier in the simulation. With practical considerations in mind, there will be some point at which costs (ethical and/or material) outweigh the benefit of testing further monkeys.

The proof of Theorem \ref{thm:eps} required a result from constructive mathematics. We conjecture that such a result is classically impossible, since the singular covering theorem is classically not true. 

A deeper fact here is that from the classical viewpoint, the computable reals have zero measure, and all finite texts produced by monkeys correspond to the rationals (or some other convenient subset of computable reals). The context of the results, then, would indicate that a careful constructive study of probability distributions provides \textit{a priori} motivation for repetition of simulations for accuracy (to rule out accidental pathological distributions generated by computer programs), and has potentially more to say about issues involving computer simulations.

There is the further issue of what model of constructive mathematics provides a good framework for this sort of work. Philosophically there is tension between the intuitionistic free choice-sequence approach and the computable sequence approach, and within these approaches are further complications by sensitivity of the theory to the validity (or invalidity) of the various versions of K\"onig's Lemma. It is not the aim here to go deeply into these issues, which could lead to a lengthy series of papers. In the interest of brevity, we leave such explorations for future research.

It has not escaped our attention that science and mathematics have each been considered to be ``games'' of recombining a finite set of characters (even if we do not yet know what they all are). Even if we consider only finite strings which are syntactically sound, and not contradicted by empirical evidence, our result shows that completing such a list is not necessarily even \emph{likely} to happen within \emph{any} finite time, such as a human lifespan, the duration of a civilisation, or even the age of the universe.

\paragraph{Acknowledgements:} The authors would like to acknowledge the contributions of the anonymous referees for substantial improvements to the paper. McKubre-Jordens was partially funded by Marsden Fund Fast-Start Grant UC1205.

\bibliographystyle{abbrv}

\begin{thebibliography}{99}
\bibitem{aristotle} Aristotle (350 BCE) \emph{Metaphysics}. Translation by W.D.\! Ross. \url{http://classics.mit.edu/Aristotle/metaphysics.html}. Retrieved 18 June 2010.
\bibitem{borel}\'{E}.\! Borel (1913) `M\'{e}canique Statistique et Irr\'{e}versibilit\'{e}'. \emph{J. Phys.}, 5e s\'{e}rie 3, 189--196.
\bibitem{borges} J.L.\! Borges (1939)  In \emph{The Total Library: Non-Fiction 1922-1986}. Translated by E.\! Weinberger (2000). Penguin, London, 214--216. 
\bibitem{varieties} D.S.\! Bridges \& F.\! Richman (1987) \emph{Varieties of Constructive Mathematics}. LMS Lecture Notes Series, Cambridge University Press.
\bibitem{chan} Y.K.\! Chan (1974) Notes on Constructive Probability Theory. \emph{Ann. Prob.}, 2(1), 51--75. Institute of Mathematical Statistics.
\bibitem{eddington} A.\! Eddington (1928) \emph{The Nature of the Physical World: The Gifford Lectures}. New York: Macmillan.
\bibitem{complete} Elmo, Gum, Heather, Holly, Mistletoe, \& Rowan \emph{Notes Towards The Complete Works of William Shakespeare} (2002) Khave-Society \& Liquid Press, UK.
\bibitem{kk} C.\! Kittel \& H.\! Kroemer (1980) \emph{Thermal Physics} (2nd ed.). W.H.\! Freeman Company.
\bibitem{kushner} B.A.\! Kushner (1999) Markov's constructive analysis; a participant's view. \emph{Theoretical Computer Science} 219, 267--285.
%\bibitem{nasa} NASA `The Beginnings of Research in Space Biology at the Air Force Missile Development Center, 1946--1952', \emph{History of Research in Space Biology and Biodynamics}. \url{http://history.nasa.gov/afspbio/part1.htm}. Retrieved 18 June 2010.
\bibitem{paignton} `Give six monkeys a computer, and what do you get? Certainly not the Bard', https://www.theguardian.com/uk/2003/may/09/science.arts.
\bibitem{shakespeare} W.\! Shakespeare \emph{The Complete Works of William Shakespeare} (2001) Geddes \& Grosset, Scotland.
\bibitem{zaslavsky} I.D.\! Zaslavsky \& G.S.\! Tseitin (1962) Singular coverings and properties of constructive functions connected with them, \emph{Problems of the constructive direction in mathematics. Part 2. Constructive mathematical analysis}, Collection of articles, Trudy Mat. Inst. Steklov., 67, Acad. Sci. USSR, Moscow–Leningrad, 458--502. English translation: A.M.S. Translations (2) 98 (1971), 41-89, MR 27\#2408.
\end{thebibliography}

\end{document}